\magnification = \magstephalf
\hsize=14truecm
\hoffset=1truecm
\parskip=5pt
\nopagenumbers 
\overfullrule=0pt
\font\eightrm=cmr8 \font\eighti=cmti8 \font\eightbf=cmbx8
\headline={\ifnum\pageno=1 \hfill\else%
{\tenrm\ifodd\pageno\rightheadline \else
\leftheadline\fi}\fi}
\def\rightheadline{EJDE--1999/11\hfil Navier-Stokes equation \hfil\folio}
\def\leftheadline{\folio\hfil Stephen Montgomery-Smith \hfil EJDE--1999/11}
\voffset=2\baselineskip
\vbox {\eightrm\noindent\baselineskip 9pt %
 Electronic Journal of Differential Equations,
Vol. {\eightbf 1999}(1999) No.~11, pp. 1--19.\hfill\break
ISSN: 1072-6691. URL: http://ejde.math.swt.edu or http://ejde.math.unt.edu
\hfil\break ftp  ejde.math.swt.edu \quad ftp ejde.math.unt.edu (login: ftp)}
\footnote{}{\vbox{\hsize=10cm\eightrm\noindent\baselineskip 9pt %
1991 {\eighti Subject Classification:} 35Q30, 76D05, 35B65.
\hfil\break
{\eighti Key words and phrases:} Navier-Stokes equation, thin domain.
\hfil\break
\copyright 1999 Southwest Texas State University  and
University of North Texas.\hfil\break
Submitted November 19, 1998. Published April 14, 1999.\hfil\break
Partially supported by grants from the N.S.F.\
and the Research Board of the University of Missouri.} }

\bigskip\bigskip

\centerline{GLOBAL REGULARITY OF THE NAVIER-STOKES EQUATION }
\centerline{ON THIN THREE-DIMENSIONAL DOMAINS }
\centerline{WITH PERIODIC BOUNDARY CONDITIONS }
\medskip
\centerline{Stephen Montgomery-Smith}
\bigskip\bigskip

{\eightrm\baselineskip=10pt \narrower
\centerline{\eightbf Abstract}
This paper gives another version of results due to Raugel and Sell, and 
similar results due to Moise, Temam and Ziane, that state the following: 
the solution of the Navier-Stokes equation on a thin three-dimensional 
domain with periodic boundary conditions has global regularity, as long as
there is some control on the size of the initial data and the forcing
term, where the control is larger than that obtainable via ``small
data'' estimates.  The approach taken is to consider the three-dimensional
equation as a perturbation of the equation when the vector field does
not depend upon the coordinate in the thin direction.
\bigskip}

\def\modo#1{{\left|#1\right|}}
\def\normo#1{{\mathopen\|#1\mathclose\|}}

\bigbreak
\centerline{\bf \S 1. Introduction} \medskip\nobreak
The celebrated Navier-Stokes equation is concerned with the
velocity vector $u$ on a domain $\Omega$, describing the flow
of an incompressible fluid.  A famous unsolved problem is
the following: if $\Omega$ is a nice enough
3-dimensional domain, and if
the initial data is smooth, and the forcing term is uniformly 
smooth in time, then does it follow that the solution is
smooth for all time?  What is known is that a weak solution
exists, although it is not known if that solution is unique.
For details, we refer the reader to a number of books and
papers, including [CF], [DG] and [T].

For the 2-dimensional problem, the solution is known.  A
heuristic approach to solving the 3-dimensional problem is as
follows: if the solution becomes less smooth, then since we
are dealing with an incompressible fluid, the complicated
activity is going to get squashed into flat sheets, and one might
hope that the solutions on the flat sheets can be somehow dealt
with by the 2-dimensional case.  Certainly, this ``flattening''
is observed in numerical and physical experiments.

For this reason, it would seem that in order to get some
handle on the real problem, that it might be important to understand
what happens to the solution to the Navier-Stokes equation on
thin domains, that is, domains of the form $M \times [0,\epsilon]$,
where $M$ is some 2-dimensional manifold, and $\epsilon>0$
is a small number.  This is what Raugel
and Sell did in a series of papers [RS1], [RS2], [RS3], as did
Avrin [A], Temam and Ziane [TZ1], [TZ2], 
Moise, Temam and Ziane [MTZ], and Iftimie [I1], [I2].

In this paper, we consider the same situation considered by
Raugel and Sell in [RS2], or by Moise, Temam and Ziane in [MTZ],
or by Iftimie in [I1], [I2].  
Let $\Omega_\epsilon = [0,l_1] \times
[0,l_2] \times [0,\epsilon]$, where $l_1 \ge l_2$ are some
positive numbers, and $\epsilon \in (0,l_2/4)$ is
some small number.
Let us consider vector fields
$u : \Omega_\epsilon
\to {\bf R}^3$ satisfying the periodic boundary
conditions
$$ u(0,y,z) = u(l_1,y,z), \quad
   u(x,0,z) = u(x,l_2,z), \quad
   u(x,y,0) = u(x,y,\epsilon) .\eqno(1)$$
Given any appropriately smooth (where being in $L_2$ is smooth
enough)
vector field $u$ satisfying these boundary conditions, 
we may split it into its divergence-free part
$Lu$ and its gradient part.  Thus $L$ is the 
so called Leray projection.

The Navier-Stokes equation considered in this paper
is the equation for a function
$u(t)$, $t\ge0$, taking values in the space of 3-dimensional
divergence-free vector fields on $\Omega_\epsilon$ satisfying
the boundary conditions (1).  The equation is
$$ \partial_t u = \nu\Delta u - L(u\cdot \nabla u) + L(f) ,\eqno(2)$$
where $\nu$ is a positive constant (the viscosity), and $f$ is a function
of $t$ taking values in the 3-dimensional vector fields satisfying
the boundary conditions (the forcing term).
For simplicity let us suppose $f = L(f)$.  

In fact, to simplify our presentation, it will make sense to
assume that \hfil\break $\int_{\Omega_\epsilon} f \, dV = 0$, and that
$\int_{\Omega_\epsilon} u(0) \, dV = 0$.  In that case it
is not hard to see that we have that $\int_{\Omega_\epsilon} u(t) \, dV = 0$
for all $t \ge 0$.  This assumption does not really affect the
solution very much, for suppose that we do not have this assumption.
Given any function $g$ on $\Omega_\epsilon$, let $\bar g$ denote its
mean value $(l_1 l_2 \epsilon)^{-1} \int_{\Omega_\epsilon} g \, dV$.
Let $(\xi_t,\eta_t,\zeta_t) = \int_0^t \bar u(0) + \bar f(s) \, ds$.
Then 
replacing $u(x,y,z,t)$ by $u(x+\xi_t,y+\eta_t,z+\zeta_t,t) - 
\bar u$ gives us
a solution to the Navier-Stokes equation in which $f$ is replaced
by $f - \bar f$, and in which $\bar u = 0$.

It is known that in order to show global regularity of
$u$ that it is sufficient to show that $u$ stays within
the Sobolev space $H^1$, that is, the space in which one
derivative is in $L_2$.  Furthermore, once this is
established, it also follows that the solution is unique.
(See [CF] or [DG].)  We will also include 
results concerning the Sobolev space $H^2$, the space in which
two derivatives are in $L_2$.

Throughout this paper, we
will use the letter $c$ to denote a positive constant (typically
larger than one), whose value will change with each occurrence.
Only in certain places (such as in Lemmas~3 and~5) will we 
use subscripts on the $c$'s to identify them, so as to avoid confusion.

\def\middlethmone{Let 
$$ U = \normo{u(0)}_{H^1}, \quad
   F = \sup_t \normo{f(t)}_2, \quad
   M = \max\left\{U,{l_1\over\nu}\,F\right\} .$$}

\def\endthmone{there exists a solution with the following
properties.  First, $u(t)$ is
in $H^1$ for all $t\ge 0$, with
$$ \normo{u(t)}_{H^1} \le c 
   \max\left\{M,
   {l_1^{3/2} \over \nu l_2^{1/2}} \, \epsilon^{-1/2} M^2\right\}. $$
If $\displaystyle t \ge c \, {l_1^2 \over \nu}$, then
$$ \normo{u(t)}_{H^1} \le 
   c 
   \max\left\{{l_1 \over \nu} F,
   {l_1^{7/2} \over \nu^3 l_2^{1/2} } \, \epsilon^{-1/2} F^2 \right\}. $$
Second, $u(t)$ is in $H^2$ for almost all $t \ge 0$, and indeed
$$ \int_0^t \normo{u(s)}_{H^2}^2 \, ds < \infty $$
for all $0\le t < \infty$.
}

\proclaim Theorem 1.  Let
$u$ satisfy the Navier-Stokes equation (2) with periodic
boundary conditions (1), and $\int_{\Omega_\epsilon} u \, dV = 0$.
\middlethmone
If $\displaystyle M \le c^{-1} {\nu l_2^{1/2} \over l_1}$, then 
\endthmone

Using a rescaling argument,
it may be shown that is is sufficient to show Theorem~1 in the case
that $l_1$, $l_2$ and $\nu$ are of order one, say, 
that these quantities all lie between $1/2$ and $2$.  We will demonstrate
this at the end of the proof of Theorem~1.  So for the remainder of this 
discussion, let us focus on this case.  Then in effect, $M = \max\{U,F\}$, 
and Theorem~1 gives global regularity in the case that $M \le c^{-1}$.

This result is not obtainable by the usual ``small data'' arguments, because
these would only give regularity in the case that $M \le c^{-1} \epsilon^{1/2}$.
(The $\epsilon^{1/2}$ comes from the fact that we are calculating
an $L_2$ norm on a domain whose total measure is of order $\epsilon$.)

We will obtain Theorem~1 by considering it as a perturbation of the
Navier-Stokes equation in which neither $u$ nor $f$ depend upon the
third coordinate of $\Omega_\epsilon$ (the `$z$' coordinate).
This approach was also taken by Iftimie [I1], [I2], and
by Moise, Temam and Ziane [MTZ].
Note that the following result effectively does not depend upon 
$\epsilon$, in that if the result is obtained for one value of $\epsilon$,
it is then automatically obtained for all others by a rescaling argument.

\proclaim Theorem 2.  Let
$u$ satisfy the Navier-Stokes equation (2) with periodic
boundary conditions (1), and $\int_{\Omega_\epsilon} u \, dV = 0$,
and neither $u$ nor $f$ depend upon the third coordinate of
$\Omega_\epsilon$.
\middlethmone
Then 
\endthmone

In order to compare Theorem~1 with the results in the literature,
let us define the following projections.
Let 
$$ Pu(x,y,z) = {1\over\epsilon} \int_0^\epsilon u(x,y,\zeta) \, d\zeta ,$$
and let $Qu = u - Pu$.
As we stated above, we will discuss these results only in the
case when $l_1$, $l_2$ and $\nu$ are of order one.
Then the result of Raugel and Sell [RS2] gives global $H^1$ boundedness
of the solution as long as
$$ \eqalignno{
   \normo{Pu(0)}_{H^1} 
   &\le \epsilon^{{7\over 24} + \delta_1} 
        (\log(1/\epsilon))^{\delta_2} \cr
   \normo{Qu(0)}_{H^1} 
   &\le \epsilon^{-{5\over 48} + \delta_3} 
        (\log(1/\epsilon))^{\delta_4} \cr
   \sup_t \normo{Pf(t)}_2 
   &\le \epsilon^{{7\over 24} + \delta_5} 
        (\log(1/\epsilon))^{\delta_6} \cr
   \sup_t \normo{Qf(t)}_2 
   &\le \epsilon^{-{1\over 2} + \delta_7} 
        (\log(1/\epsilon))^{\delta_8} ,\cr } $$
where $\delta_i$ ($1\le i \le 8$) are positive numbers.

The result of Moise,
Temam and Ziane [MTZ] gives $H^1$ boundedness for $t \in [0,T(\epsilon)]$,
where $\lim_{\epsilon\to0}T(\epsilon) = \infty$, and also
an integral condition on the $H^2$ norm, as long as
$$ \eqalignno{
   \normo{Pu(0)}_{H^1} 
   &\le \alpha(\epsilon)\epsilon^{{1\over 6} + \delta} 
         \cr
   \normo{Qu(0)}_{H^1} 
   &\le \alpha(\epsilon)\epsilon^{-{1\over 6} + \delta} 
        \cr
   \sup_t \normo{Pf(t)}_2 
   &\le \alpha(\epsilon)\epsilon^{{1\over 6} + \delta} 
        \cr
   \sup_t \normo{Qf(t)}_2 
   &\le \alpha(\epsilon)\epsilon^{-{1\over 6} + \delta} 
        \cr } $$
where $\delta$ is a positive number, and 
$\lim_{\epsilon\to0}\alpha(\epsilon)=0$.

Iftimie [I1], [I2] gets global existence results 
under conditions that use certain `aniso\-tropic' Sobolev spaces.
For example, one case of his Theorem~4.1 gives
global existence under the condition that the forcing term $f$ is zero, and
$$ \normo{Qu}_{H^{1/2}} \exp(c \epsilon^{-1} \normo{Pu}_2^2) \le c^{-1} .$$
Even though his conclusions are slightly different, it is instructive to see
how his hypothesis relates to that of this paper.  
Indeed, his condition is true if
$$ \eqalignno{
   \normo{Pu}_{H^1} 
   &\le c_\delta^{-1} \epsilon^{1/2} \sqrt{\log(1/\epsilon)} 
        \cr
   \normo{Qu}_{H^1} 
   &\le c^{-1} \epsilon^{-1/2 + \delta} , \cr } $$
where $c_\delta$ depends upon $\delta>0$.

Before proceeding further, let us set our notation, and also 
state the ``tools of the trade,'' that is, the standard results
that many people in this area use.  
(See [CF] or [DG] --- in particular [DG] considers the case of
periodic boundary conditions, and many of the following
calculations may be found there.)

As well as the operators $P$ and $Q$ given
above, let us define the operators $R$ and $S$:
$R((u_1,u_2,u_3)) = (u_1,u_2,0)$, and $Su = u - Ru = (0,0,u_3)$.
Split $u = v + w = r + s + w$, where $v = Pu$, $w = Qu$, $r = Rv$, 
$s = Sv$.  Since
$r$ and $s$ do not depend upon $z$ (where we label the coordinates
of $\Omega_\epsilon$ by $x$, $y$ and $z$), it is clear that
$r$, $s$, $v$ and $w$ are all divergent free vector fields.

As a notational device, I will write
$$ \normo{D^\alpha u}_p $$
for the Sobolev space $\alpha$ derivatives in $L_p$.  (Thus $D$ might
represent the operator $\sqrt{-\Delta}$.)
We have the Sobolev inequalities: if
$f$ is a mean zero
function on $[0,l_1]\times [0,l_2]$, $1 < p < q < \infty$, and
$\alpha > 0$, then
$$ \normo{f}_q \le c \normo{D^\alpha f}_p ,$$
where ${1\over q} + {\alpha \over 2} = {1\over p}$, and
$c$ depends upon $p$ and $q$, as well as $l_1$ and $l_2$.  Thus, if
$f$ is a mean zero function on $\Omega_\epsilon$, then
$$ \normo{Pf}_q \le c \epsilon^{-\alpha/2} \normo{D^\alpha Pf}_p .$$
In fact, the only condition under which we will use this inequality is
in the case $p=4$ and $q=2$, when $\alpha = 1/2$.  For this case,
we will include an elementary proof in the Appendix.

We have the interpolation inequality: if $f$ is a mean zero
function,
$\alpha_0,\alpha_1$ are real numbers, and $0 \le \theta \le 1$, then
$$ \normo{D^{\alpha_\theta} f}_2 \le c
   \normo{D^{\alpha_0} f}_2^{1-\theta}
   \normo{D^{\alpha_1} f}_2^\theta ,$$
where $\alpha_\theta = (1-\theta)\alpha_0 + \theta \alpha_1$.
This inequality is easy to
show using Parseval's identity and H\"older's inequality.
(See for example the proofs of Lemmas~4
or~6 for the statement of Parseval's identity.)

We also have the Poincar\'e inequalities: if $f$ is a mean zero function
on $\Omega_\epsilon$,
with periodic boundary conditions,
and $\alpha>0$, then
$$ \normo{f}_2 \le c \normo{D^\alpha f}_2 ,$$
where $c$ depends upon $\alpha$ as well as $l_1$ and $l_2$.
Again, 
this is easy to show using Parseval's identity.

If $u$ is a divergence-free vector field on the domain 
$\Omega_\epsilon$
with periodic boundary conditions, and if $f$ and $g$ are two other 
functions on $\Omega_\epsilon$ with periodic boundary conditions, 
sufficiently smooth so that the
following integrals make sense, then by integration by parts we get
$$ \int_{\Omega_\epsilon} f (u \cdot \nabla g) \, dV
   =
   - \int_{\Omega_\epsilon} g (u \cdot \nabla f) \, dV ,$$
and so
$$ \int_{\Omega_\epsilon} f (u \cdot \nabla f) \, dV = 0 .$$

If $r$ is a two-dimensional, divergence-free vector field on the domain
$[0,l_1] \times [0,l_2]$ with periodic boundary conditions that
is sufficiently smooth, then we have the ``enstrophy miracle:''
$$ \int_{[0,l_1] \times [0,l_2]} 
   \Delta r \cdot (r \cdot \nabla r) \, dA = 0 .$$
This is obtained as follows.  First, integrating by parts, we see that
the left hand side is equal to
$$ \eqalignno{
   &- \int_{[0,l_1] \times [0,l_2]}
     \partial_x r \cdot (\partial_x r \cdot \nabla r) \, dA
   - \int_{[0,l_1] \times [0,l_2]} 
     \partial_x r \cdot (r \cdot \nabla \partial_x r) \, dA \cr
   &\quad- \int_{[0,l_1] \times [0,l_2]} 
     \partial_y r \cdot (\partial_y r \cdot \nabla r) \, dA
   - \int_{[0,l_1] \times [0,l_2]} 
     \partial_y r \cdot (r \cdot \nabla \partial_y r) \, dA .\cr}$$
(Here, as in the rest of the paper, $\partial_x$, $\partial_y$ and
$\partial_z$ represent partial differentiation with respect to
$x$, $y$ and $z$ respectively, that is, the first, second and third
coordinates respectively.)
We see that the second and fourth terms are zero.  Expanding
and collecting the first and third
terms, and remembering that $r$ is divergence-free, we see that they also
total to zero.

\beginsection Navier-Stokes for Flows Independent of $z$

Let us start with the proof of Theorem~2,
the case when $w=Qu=0$, that is, when $u = Pu = v$.  
We will prove Theorem~2 in the case that the quantities $l_1$, $l_2$ and
$\nu$ all lie between $1/2$ and $2$.  The general result may be
obtained as shown at the end of the proof of Theorem~1.

First we need to find a solution to the Navier-Stokes equation.
This is done using so called Galerkin solutions.  Let $S_n$ denote
the projection that takes a function $f$ on $\Omega_\epsilon$ onto
the $n$th partial Fourier series.  (Quite how we index this sum is
not important, as long as $S_n$ converges formally to the identity.)
Then we consider the solution $u_n$ to the problem
$$ \partial_t u_n = \nu\Delta u_n - S_n L(u_n\cdot \nabla u_n) + S_n L(f) ,$$
with $u_n(0) = S_n u(0)$ for which $u = S_n u$.  It is a well known
argument to show that this equation, essentially an ODE, has a global
solution, and that the solutions $u_n$ converge weakly to some function
$u$.  This is the so called weak solution to the Navier-Stokes equation.
In that case, for any appropriate norm $\normo\cdot$, we will have that
$\normo u \le \liminf_{n\to\infty}\normo{u_n}$.

Thus, in order to prove our theorem, it is sufficient to prove it
for the Galerkin solutions.  This is what we shall do.  However, carrying
the symbol $S_n$ throughout the proof could be a little cumbersome,
and so for this reason, we will replace $L$ by $S_n L$, and suppose
that both $f$ and $u(0)$ lie in the range of $S_n$.  (We will also
suppose that $f$ lies in the range of $L$.)

So let us proceed.  Notice that the Navier-Stokes equation becomes
$$ \partial_t v = \nu\Delta v - L(r\cdot\nabla v) + f ,$$
since $v \cdot \nabla v = r \cdot \nabla v$, because $\partial_z v = 0$.
If we apply the $R$ and $S$ operator to this, we get the following
pair of equations:
$$ \eqalign{
   \partial_t r &= \nu\Delta r - L(r\cdot\nabla r) + Rf \cr
   \partial_t s &= \nu\Delta s - L(r\cdot\nabla s) + Sf . \cr } \eqno(3)$$
The first equation is merely the 2-dimensional Navier-Stokes for
the flow $r$.  The second equation essentially says that the
1-dimensional quantity $s$ is being pushed around by the 2-dimensional
flow $r$ (and indeed in the second equation, the operator $L$ acts
as the identity).

Let us write 
$$ \phi = \normo{Dr}_2, \quad \psi = \normo{Ds}_2,
   \quad \tilde\phi = \normo{D^2r}_2, \quad \tilde\psi = \normo{D^2s}_2,
   \quad \theta = \normo u_2 .$$
Poincar\'e's inequality tells us immediately that $\phi \le c \tilde \phi$,
$\psi \le c \tilde \psi$, and $\theta^2 \le c(\phi^2+\psi^2)$.

The process for comprehending $\phi$ and $\tilde\phi$ 
is well known.
Start with the first equation in (3),
dot product both sides with $-\Delta r$, integrate over $\Omega_\epsilon$, 
use the self-adjointness of $L$,
apply
some integration by parts, use the Cauchy-Schwartz inequality, and use the 
``enstrophy miracle,''
to get
$$ {1\over 2} \partial_t \normo{Dr}_2^2 \le -\nu\normo{D^2r}_2^2 
   + \normo{D^2r}_2 F .$$

Use the inequality $ab \le a^2 + b^2$ to get that 
$\normo{D^2r}_2 F \le (\nu/2) \normo{D^2r}_2^2 + (2/\nu) F^2$.
Thus we have that
$$ \partial_t \phi^2 \le - c^{-1} \tilde\phi^2 + c F^2 .$$
This differential inequality is easy to solve, but before we
do so, let us first understand $\psi$ and $\tilde\psi$.  
Take the second equation
from (3), dot product both sides with $-\Delta s$, integrate over
$\Omega_\epsilon$, and work as before.  But in this case,
the ``enstrophy miracle'' does
not work --- there is a term:
$$ \int_{\Omega_\epsilon} \Delta s \cdot (r \cdot \nabla s) \, dV .$$
To get a grip on this term, see that it splits into 
$$ \int_{\Omega_\epsilon} \partial_x^2 s \cdot (r \cdot \nabla s) \, dV $$
plus another term with $\partial_y^2$ in place of $\partial_x^2$.  Integrate
by parts to get
$$ - \int_{\Omega_\epsilon} \partial_x s\cdot(\partial_x r \cdot\nabla s) \, dV
   - \int_{\Omega_\epsilon} \partial_x s\cdot(r \cdot \nabla\partial_x s) \, dV
   .$$
The second term is zero.  For the first term, we may use
H\"older's inequality and the Sobolev inequality to bound it by:
$$ \eqalignno{
   \normo{Dr}_2 \normo{Ds}_4^2 
   &\le c \epsilon^{-1/2} \normo{Dr}_2 \normo{D^{3/2}s}_2^2 \cr
   &\le c \epsilon^{-1/2} \normo{Dr}_2 \normo{Ds}_2 \normo{D^2s}_2 \cr
   &\le {1\over 2}\normo{D^2s}_2^2 
       + c \epsilon^{-1} \normo{Dr}_2^2 \normo{Ds}_2^2 ,\cr}$$
where in the last step we use the inequality $ab \le a^2 + b^2$.
   
Putting this all together, we get a differential inequality:
$$ \partial_t \psi^2 \le - c^{-1} \tilde\psi^2 
   + c \epsilon^{-1} \phi^2 \psi^2 + c F^2 .$$

We will also require a differential equation for $\theta$:
take the Navier-Stokes equation, dot product both sides with $u$,
integrate over $\Omega_\epsilon$, and do the usual stuff,
to get
$$ {1\over 2} \partial_t \normo u_2^2 \le -\nu\normo{Du}_2^2 
                                          + \normo{Du}_2 F .$$
Since $\normo{Du}_2 F \le (\nu/2) \normo{Du}_2^2 + (2/\nu) F^2$, we get
$$ {1\over 2} \partial_t \normo u_2^2 \le -{\nu\over 2}\normo{Du}_2^2 +
                                           {2\over\nu} F^2 ,$$
that is
$$ \partial_t \theta^2 \le -c^{-1} (\phi^2+\psi^2) + c F^2 .$$

Then Theorem~2 will be established once we have obtained the following
result.

\def\protolemma#1#2#3#4#5#6#7{%  
#7
Let $\phi$, $\psi$, $\tilde\phi$, $\tilde\psi$,
$\theta$ be positive differentiable functions
of $t$.
Suppose that
for some positive constants $c_i$ ($1\le i \le 10$#4) we have
$$ \eqalignno{
   \phi(0) & \le U & (#1.1) \cr
   \psi(0) & \le U & (#1.2) \cr
   \theta^2 & \le c_1 (\phi^2 + \psi^2) & (#1.3)\cr
   \phi & \le c_2 \tilde \phi & (#1.4)\cr
   \psi & \le c_3 \tilde \psi & (#1.5)\cr
   \partial_t \phi^2 &\le 
      #2
      - c_4^{-1} \tilde \phi^2 + c_5 F^2 & (#1.6)\cr
   \partial_t \psi^2 &\le 
      #2
      - c_6^{-1} \tilde \psi^2 
      + c_7 \epsilon^{-1} \phi^2 \psi^2 + c_8 F^2 & (#1.7)\cr
   \partial_t \theta^2 &\le -c_9^{-1} (\phi^2+\psi^2) + c_{10} F^2 
           & (#1.8)\cr}$$
for $0 \le t < #6$.
Let $M = \max\{U,F\}$.
Then there exist positive constants $c_i$ ($11\le i \le 17$#5), depending
only upon $c_i$ ($1 \le i \le 10$) such that 
#3
we have the inequalities
$$ \eqalignno{
   \theta^2 &\le c_{11}(F^2 + (U^2 - F^2) e^{-c_{12}^{-1}t}) & (#1.9)\cr
   \phi^2 &\le c_{13}(F^2 + (U^2 - F^2) e^{-c_{14}^{-1}t}) & (#1.10)\cr
   \psi &\le c_{15} \max\{\epsilon^{-1/2} M^2,M\}, & (#1.11)\cr
\noalign{\noindent and if $t \ge c_{16}$ then}
   \psi &\le c_{17} 
      \max\{\epsilon^{-1/2} F^2, F \}. & (#1.12)\cr}$$
for $0 \le t < #6$.  Furthermore, we have that
$$ \int_0^t (\phi(s)^2+\psi(s)^2) \, ds < \infty \eqno (#1.13) $$
for $0 \le t < #6$.
}

\proclaim Lemma 3. \protolemma{4}{}{}{}{}{T}
{Let $U$, $F$ and $T$ be positive numbers, where $T$ may be infinity.}

\noindent Proof:\ \ 
Inequalities (4.9) and (4.10) are easy to obtain from
combining (4.1) and (4.4) with (4.6), and (4.1), (4.2) and (4.3) with (4.8),
by using
Gronwall's
inequality.  

Let us obtain (4.11).
 From (4.10), we see that there is a positive number 
$c_{18}$ such that
$\phi \le c_{18} M$.
Combining this with inequalities (4.5), (4.7) and (4.8), we see that
for some positive constants $c_{19}$ and $c_{20}$ that
$$ \partial_t \psi^2 + c_{19}\psi^2 
   \le c_{20} \epsilon^{-1} M^2 (F^2 - \partial_t \theta^2)
   + c_{20} F^2 .$$
Multiply both sides by $e^{c_{19} t}$ and integrate from $0$ to $t$ to
get that
$$ 
   e^{c_{19} t} \psi^2 - \psi(0)^2
   \le c_{20}
   \int_0^t \epsilon^{-1} M^2 e^{c_{19}s} (F^2 - \partial_s \theta(s)^2)
            + F^2 e^{c_{19}s} \, ds ,$$
which, by evaluating the integrals, and integrating by parts, is less than
or equal to
$$  c_{20} c_{19}^{-1} e^{c_{19}t} (\epsilon^{-1} M^2 F^2 + F^2)
     + c_{20} \epsilon^{-1} M^2 (\theta(0)^2 - e^{c_{19}t}\theta^2)
     + c_{20} c_{19} \epsilon^{-1} M^2 
       \int_0^t e^{c_{19}s} \theta(s)^2 \, ds. $$
Now, from (4.9), we see that there is a positive constant $c_{21}$ such that
$\theta \le c_{21} M$.  Hence
$$ \eqalignno{
   e^{c_{19} t} \psi^2 - \psi(0)^2
   & \le c_{20} c_{19}^{-1} e^{c_{19}t} (\epsilon^{-1} M^2 F^2 + F^2)
     + c_{20} c_{21}^2 \epsilon^{-1} M^4
     + c_{20} c_{19} \epsilon^{-1} M^2 c_{21}^2 \int_0^t e^{c_{19}s} \, ds \cr
   & \le c_{20} c_{19}^{-1} e^{c_{19}t} (\epsilon^{-1} M^2 F^2 + F^2)
     + c_{20} c_{21}^2 \epsilon^{-1} M^4
     + c_{20} \epsilon^{-1} M^2 c_{21}^2 e^{c_{19}t} . \cr } $$
Hence
$$ \psi^2 \le e^{-c_{19} t} U 
   + c_{20} c_{19}^{-1}(\epsilon^{-1} M^2 F^2 + F^2)
   + c_{20} c_{21}^2 \epsilon^{-1} e^{-c_{19} t} M^4
     + c_{20} c_{21}^2 \epsilon^{-1} M^2 ,$$
and from here it is easy to obtain (4.11).

To obtain (4.12) is similar.  We see that there are positive numbers
$c_{22}$, $c_{23}$ and $c_{24}$ such that if $t \ge c_{22}$ then
$\phi \le c_{23} F$ and $\theta \le c_{24} F$.  Apply the above argument,
except integrate from $c_{22}$ to $t$ instead of from $0$ to $t$.

Finally, (4.13) may be obtained by integrating (4.6) and (4.7), and using
(4.9), (4.10) and (4.11).
\hfill Q.E.D.

\beginsection Proof of Theorem~1

The proof of Theorem~1 will follow the same lines as the proof of Theorem~2,
with some additional work for dealing with the $w = Qu$ part.  
Let us start by assuming that the quantities $l_1$, $l_2$ and $\nu$
all lie between $1/2$ and $2$.

We need a couple of Poincar\'e/Sobolev type inequalities on
$\Omega_\epsilon$.

\proclaim Lemma 4.  Let $w=Qu$, then
$$ \eqalignno{
   \normo w_\infty &\le c \epsilon^{1/2} \normo{D^2 w}_2 \cr
   \normo w_4 &\le c \epsilon^{1/4} \normo{D w}_2 .\cr } $$

\noindent Proof:\ \ 
Let $\hat w$ denote the Fourier coefficients of $w$:
$$ \hat w(m,n,p) = (l_1 l_2 \epsilon)^{-1} \int_{\Omega_\epsilon}
   w(x,y,z) \exp(-2\pi i(mx/l_1+ny/l_2+pz/\epsilon) \, dx\,dy\,dz ,$$
where $m$, $n$ and $p$ are integers.  
Then the function $w$ can be reconstructed using the Fourier series
$$ w(x,y,z) = \sum_{m,n,p} \hat w(m,n,p)
  \exp(2\pi i(mx/l_1+ny/l_2+pz/\epsilon) .$$
We recall Parseval's identity:
$$ \normo w_2^2 = l_1 l_2 \epsilon
   \sum_{m,n,p} \modo{\hat w(m,n,p)}^2 ,$$
and the Hausdorff-Young inequality: if $2 \le p \le \infty$, and 
$p' = p/(p-1)$, then
$$ \normo w_p \le (l_1 l_2 \epsilon)^{1/p}
   \left( \sum_{m,n,p} \modo{\hat w(m,n,p)}^{p'}\right)^{1/p'} .$$
Since $w = Qu$, we have that $\hat w(m,n,p) = 0$ if $p = 0$.

Now for any real number $\alpha$, we have that
$$ \widehat{D^\alpha w}(m,n,p) = (-2\pi i)^\alpha
   (m^2/l_1^2+n^2/l_2^2+p^2/\epsilon^2)^{\alpha/2}
   \hat w(m,n,p) ,$$
and thus by Parseval's identity
$$ \normo{D^\alpha w}_2 = (2 \pi)^\alpha (l_1 l_2 \epsilon)^{1/2} 
   \left( \sum_{m,n,p} (m^2/l_1^2+n^2/l_2^2+p^2/\epsilon^2)^\alpha
          \modo{\hat w(m,n,p)}^2 \right)^{1/2} . $$

Let us start with showing the first inequality.  Apply Cauchy-Schwartz
to get
$$ \eqalignno{
   &\normo w_\infty 
   \le 
   \sum_{m,n,p} \modo{\hat w(m,n,p)} \cr
   &\quad\le
   \left(\sum_{m,n,p} 
   I_{p\ne 0} (m^2/l_1^2+n^2/l_2^2+p^2/\epsilon^2)^{-2} \right)^{1/2} \times\cr
   &\qquad\quad
   \left(\sum_{m,n,p} (m^2/l_1^2+n^2/l_2^2+p^2/\epsilon^2)^2
          \modo{\hat w(m,n,p)}^2 \right)^{1/2} .\cr}$$
By approximating sums by integrals, and using other elementary
inequalities, we see that
$$ \eqalignno{
   \sum_{m,n,p} I_{p\ne 0} (m^2/l_1^2+n^2/l_2^2+p^2/\epsilon^2)^{-2} 
   &\le c
   \sum_{m,n,p} I_{p\ne 0} {1\over\max\{m^4,n^4,p^4/\epsilon^4\}} \cr
   &\le c
   \sum_{n,p} I_{p \ne 0} {1\over\max\{n^3,p^3/\epsilon^3\}} \cr
   & \le c
   \sum_p I_{p \ne 0} {1\over p^2/\epsilon^2} \cr
   & \le c \epsilon^2 .\cr }$$
Hence we obtain the first inequality.

The second inequality has a similar proof: start by
using H\"older's inequality to get
$$ \eqalignno{
   &\normo w_4 
   \le (l_1 l_2 \epsilon)^{1/4}
   \left( \sum_{m,n,p} \modo{\hat w(m,n,p)}^{4/3} \right)^{3/4} \cr
   &\quad\le (l_1 l_2 \epsilon)^{1/4}
   \left(\sum_{m,n,p} 
   I_{p\ne 0} (m^2/l_1^2+n^2/l_2^2+p^2/\epsilon^2)^{-2} \right)^{1/4} \times\cr
   &\qquad\quad
   \left(\sum_{m,n,p} (m^2/l_1^2+n^2/l_2^2+p^2/\epsilon^2)
          \modo{\hat w(m,n,p)}^2 \right)^{1/2} ,\cr}$$
and proceed as with the proof of the first inequality.
\hfill Q.E.D.

\bigskip

\noindent Proof of Theorem~1:\ \ 
As in the proof of Theorem~2, we argue that we work with the
Galerkin approximations.
We will obtain differential inequalities for the following
quantities:
$$ \phi = \sqrt{\normo{Dr}_2^2 + \normo{Dw}_2^2}, 
   \quad \psi = \sqrt{\normo{Ds}_2^2 + \normo{Dw}_2^2},$$
$$ \tilde\phi = \sqrt{\normo{D^2r}_2^2+\normo{D^2w}_2^2}, 
   \quad \tilde\psi = \sqrt{\normo{D^2s}_2^2+\normo{D^2w}_2^2},$$
$$ \chi = \normo{D^2w}_2,
   \quad \theta = \normo u_2 .$$
Poincar\'e's inequality tells us immediately that $\phi \le c \tilde \phi$,
$\psi \le c \tilde \psi$, and $\theta^2 \le c(\phi^2+\psi^2)$.

Let us start with the Navier-Stokes equation, and apply the operator
$P$.  Note that if $f$ and $g$ are functions on $\Omega_\epsilon$, then
$P((Pf)(Qg)) = 0$.  Thus, we obtain
 
$$ \partial_t v = \nu\Delta v - LP(v\cdot\nabla v) - LP(w\cdot\nabla w) + Pf
   = \nu\Delta v - LP(r\cdot\nabla v) - LP(w\cdot\nabla w) + Pf .$$
Now apply $R$ to both sides:
$$ \partial_t r = \nu\Delta r 
                  - LP(r\cdot\nabla r) -LP(w\cdot\nabla Rw) + PRf. $$
Finally, take the dot product of both sides with $-\Delta r$, and integrate
over $\Omega_\epsilon$, and do all the usual stuff.  A lot of the terms
work in exactly the same way that they did in the previous section.  The only
term that we did not deal with is the following:
$$ \int_{\Omega_\epsilon} \Delta r \cdot (w\cdot\nabla Rw) \, dV .$$
This splits into a term:
$$ \int_{\Omega_\epsilon} \partial_x^2 r \cdot (w\cdot\nabla Rw) \, dV ,$$
and a similar one with $\partial_y^2$ in place of $\partial_x^2$.  
The bounds for the second term will be as for the first, so let us only deal
with the first.
Integrate
by parts to get
$$ -\int_{\Omega_\epsilon} \partial_x r\cdot(\partial_x w\cdot \nabla Rw) \,dV
   -\int_{\Omega_\epsilon} \partial_x r\cdot(w\cdot \nabla\partial_x Rw) \,dV
   .$$
The first term is bounded by
$$ c \normo{Dv}_2 \normo{Dw}_4^2 ,$$
and the second by
$$ c \normo{Dv}_2 \normo{D^2w}_2 \normo w_\infty,$$
Combining all this with Lemma~4, we get
$$ \partial_t \normo{D r}_2^2
   \le
   - c^{-1} \normo{D^2 r}_2^2 + c \epsilon^{1/2} \normo{Dv}_2 \normo{D^2 w}_2^2 
   + c F^2 .
    \eqno(5)$$

The work for $s$ is practically identical, and we get
$$ \partial_t \normo{D s}_2^2
   \le
   - c^{-1} \normo{D^2 s}_2^2 
   + c \epsilon^{-1} \normo{Dr}_2^2 \normo{D s}_2^2
   + c \epsilon^{1/2} \normo{Dv}_2 \normo{D^2 w}_2^2 +c F^2 .\eqno(6)$$

We also need to establish an equation for $w$.  Take the Navier-Stokes
equation and apply $Q$.  Note that if $f$ and $g$ are two functions on
$\Omega_\epsilon$, then $Q((Pf)(Pg)) = 0$.  Thus we get
$$ \partial_t w = \nu\Delta w - LQ(w\cdot \nabla v) - LQ(v \cdot \nabla w)
                  - LQ(w \cdot \nabla w) + Qf .$$
Take the dot product with $-\Delta w$, and integrate over $\Omega_\epsilon$,
doing all the stuff as before.  Let us see what happens to the non-linear
terms, only bothering with the $\partial_x^2 w$ part of $\Delta w$, knowing
that the other parts will give the same estimates.

First we get
$$ \int_{\Omega_\epsilon} \partial_x^2 w \cdot (w\cdot \nabla v) \, dV
   =
   - \int_{\Omega_\epsilon} \partial_x w \cdot(\partial_x w\cdot\nabla v) \, dV
   - \int_{\Omega_\epsilon} \partial_x w \cdot(w\cdot \nabla\partial_x v) \, dV
   .$$
The first term is bounded in absolute value by
$$ \normo{Dw}_4^2 \normo{Dv}_2 \le c \epsilon^{1/2} \normo{D^2w}_2^2
   \normo{Dv}_2 .$$
The second term is equal to
$$ \int_{\Omega_\epsilon} \partial_x v \cdot(w\cdot\nabla\partial_x w) \, dV
   ,$$
which is bounded in absolute value by
$$ \normo{w}_\infty \normo{D^2w}_2 \normo{Dv}_2 \le 
   c \epsilon^{1/2} \normo{D^2w}_2^2
   \normo{Dv}_2 .$$

Next, we have
$$ \int_{\Omega_\epsilon} \partial_x^2 w \cdot(v\cdot \nabla w) \, dV
   =
   - \int_{\Omega_\epsilon} \partial_x w \cdot(\partial_x v\cdot \nabla w) \, dV
   - \int_{\Omega_\epsilon} \partial_x w \cdot(v\cdot \nabla\partial_x w) \, dV
   .$$
The first term is bounded in absolute value by
$$ \normo{Dw}_4^2 \normo{Dv}_2 \le c \epsilon^{1/2} \normo{D^2w}_2^2
   \normo{Dv}_2 ,$$
and the second term is zero.

Finally we have
$$ \int_{\Omega_\epsilon} \partial_x^2 w \cdot (w\cdot \nabla w) \, dV
   =
   - \int_{\Omega_\epsilon} \partial_x w \cdot(\partial_x w\cdot \nabla w) \, dV
   - \int_{\Omega_\epsilon} \partial_x w \cdot(w\cdot \nabla\partial_x w) \, dV
   .$$
The first term is bounded in absolute value by
$$ \normo{Dw}_4^2 \normo{Dw}_2 \le c \epsilon^{1/2} \normo{D^2w}_2^2
   \normo{Dw}_2 ,$$
and the second term is zero.

So, doing all the same stuff as above, we get
$$ \partial_t \normo{D w}_2^2
   \le
   - \normo{D^2 w}_2^2 
   + c \epsilon^{1/2} \normo{Dv}_2 \normo{D^2 w}_2^2 
   + c \epsilon^{1/2} \normo{Dw}_2 \normo{D^2 w}_2^2 +c F^2 .\eqno(7)$$

If we add equations (5) and (7), and also (6) and (7), 
(and apply liberally the inequalities $\sqrt{a^2 + b^2} \le a+b
\le \sqrt2 \sqrt{a^2+b^2}$ for positive $a$ and $b$)
we get the two
differential inequalities:
$$ \eqalign{
   \partial_t \phi^2 &\le -c^{-1}\tilde \phi^2 -c^{-1}\chi^2
   + c \epsilon^{1/2} (\phi+\psi) \chi^2 + c F^2\cr
   \partial_t \psi^2 &\le -c^{-1}\tilde \psi^2 -c^{-1}\chi^2
   + c \epsilon^{-1} \phi^2 \psi^2 
   + c \epsilon^{1/2} (\phi+\psi) \chi^2 + c F^2 . \cr } 
   $$
In addition, arguing as in the previous section, 
we get the differential inequality
$$ \partial_t \theta^2 \le -c^{-1}(\phi^2+\psi^2) + c F^2 $$
Thus the theorem will be established when we have proved the following
lemma.

\proclaim Lemma 5. \protolemma{8}
{(-c_{18}^{-1}+c_{19}\epsilon^{-1/2}(\phi+\psi))\chi}
{if $M \le c_{20}^{-1}$, then}
{, $18\le i \le 19$}
{, $i=20$}
{\infty}
{Let $U$ and $F$ be positive numbers.}

\noindent Proof:\ \ 
Let $(c_i)_{11\le i \le 17}$ depend upon $(c_i)_{1 \le i \le 10}$
as in Lemma~3.  Let us suppose that $U,F \le c_{20}^{-1} $,
where $c_{20}$ will be chosen momentarily.
Let
$$ T = \inf\{t > 0 : c_{19} \epsilon^{1/2} (\phi+\psi) > c_{18}^{-1}\} .$$
Suppose for a contradiction that $T < \infty$.
But then for $t \in [0,T]$, the quantities $\phi$, $\psi$, $\tilde\phi$,
$\tilde\psi$ and $\theta$
satisfy the hypothesis of Lemma~3.  But then by the conclusion of Lemma~3,
we know that for some constant $c_{21}>0$ that 
$$ \phi + \psi \le c_{21} \epsilon^{-1/2} c_{20}^{-2} .$$
Setting $c_{20}$ small enough, we see then that for $t \in [0,T]$ that
we have 
$$ c_{19} \epsilon^{1/2} (\phi+\psi) \le  c_{18}^{-1}/2 ,$$
and thus that there is a neighborhood of $T$ such that
$$ c_{19} \epsilon^{1/2} (\phi+\psi) \le  c_{18}^{-1} ,$$
contradicting the definition of $T$.

Thus $T = \infty$, and thus the functions $\phi$, $\psi$, $\tilde\phi$,
$\tilde\psi$ and $\theta$
satisfy the hypothesis of Lemma~3, and the result follows.

\bigskip

Now let us relax the restriction that $l_1$, $l_2$ and $\nu$ all lie
between $1/2$ and $2$.  Let $n$ be the integer part of 
$\displaystyle{ l_1\over l_2}$, and define
new vector fields $\tilde u$ and $\tilde f$ on 
$\displaystyle[0,1]\times\left[0,{n l_2\over l_1}\right]
\times\left[0,{\epsilon \over l_1}\right]$ according to the
formulae
$$ \eqalignno{
   \tilde u(x_1,x_2,x_3,t) &=
   {l_1 \over \nu} \ \
   u\left(l_1 x_1, (l_1 x_2 \bmod l_2), l_1 x_3, {l_1^2 \over \nu} t\right) ,\cr
   \tilde f(x_1,x_2,x_3,t) &=
   {l_1^3 \over \nu^2} \ \
   f\left(l_1 x_1, (l_1 x_2 \bmod l_2), l_1 x_3, {l_1^2 \over \nu} t\right)
.\cr}$$
Then it may be easily verified that these satisfy the equation
$$ \partial_t \tilde u = \Delta \tilde u - L(\tilde u\cdot \nabla \tilde u) 
   + L(\tilde f) ,$$
that is, one may apply the version of Theorem~1 that we already have to the
functions $\tilde u$ and $\tilde f$.  Obtaining the more general version
of Theorem~1 is then merely a question of interpreting what it says about
$\tilde u$ and $\tilde f$, taking into account the following identities:
$$ \eqalignno{
   \normo f_2 &= {\nu^2 \over n^{1/2} l_1 ^{3/2}} \normo{\tilde f}_2 , \cr
   \normo u_{H^1} &= {\nu \over n^{1/2} l_1 ^{1/2}} \normo{\tilde u}_{H^1}. 
   \cr}$$

\hfill Q.E.D.

\beginsection Appendix: the Sobolev Inequality

The following result is essentially part of the literature.  For example,
in [S]
this result is found for functions on Euclidean space.  However, for our
special case, we are able to provide an elementary proof (a
proof motivated by Littlewood-Paley theory).

\proclaim Lemma 6.  Let $f$ be a function on $[0,l_1]\times[0,l_2]$ 
satisfying periodic
boundary conditions, that is mean zero.  Then there is a positive constant $c$,
depending only upon $l_1$ and $l_2$, such that
$$ \normo f_4 \le c \normo{D^{1/2} f}_2 .$$

\noindent Proof:
For each $r = (r_1,r_2)$ a pair of integers, 
write $\modo r = \sqrt{r_1^2/l_1^2 + r_2^2/l_2^2}$.
Define the Fourier coefficients
of $f$:
$$ \hat f_r = 
   (l_1 l_2)^{-1}
   \int_0^{l_1} \int_0^{l_2} f(x,y) \exp(-2\pi i(r_1 x/l_1 + 
   r_2 y/l_2)) \, dy \, dx .$$
The original function can be reconstructed using the Fourier series
$$ f(x,y) = \sum_r \hat f_r \exp(2\pi i(r_1 x/l_1 + 
   r_2 y/l_2)) ,$$
and we have Parseval's identity
$$ \normo f_2^2 = (l_1 l_2)
   \sum_r |\hat f_r|^2 .$$
We see that
$$ \widehat{D^{1/2}f}_r = \sqrt{-2 \pi i} (r_1^2/l_1^2 + r_2^2/l_2^2)^{1/4}
   \hat f_r ,$$
and so by Parseval's equality we see that
$$ \normo{D^{1/2}f}_2^2 = 2\pi (l_1 l_2) 
   \sum_r (r_1^2/l_1^2 + r_2^2/l_2^2)^{1/2} 
   \modo{f_r}^2 .$$
For $m$ a non-negative integer, set
$$ A_m = \left(\sum_{2^m \le \modo r < 2^{m+1}} \modo{f_r}^2\right)^{1/2} .$$
Notice then that
$$ \sum_{m=0}^\infty 2^m A_m^2 \le c \normo{D^{1/2}f}_2^2  .$$
Now,
$$ \normo f_4^4 = \int_0^{l_1} \int_0^{l_2} f(x,y)^2 \overline{f(x,y)^2}
   \, dy \, dx ,$$
and expanding this using the Fourier series, we obtain that
$$ \eqalignno{
   \normo f_4^4 
   &= 
   (l_1 l_2)
   \sum_{r^{(1)}+r^{(2)}-r^{(3)}-r^{(4)} = 0}
	 \hat f_{r^{(1)}} \hat f_{r^{(2)}} 
	 \overline{\hat f_{r^{(3)}} \hat f_{r^{(4)}}} \cr
   &\le
   24 (l_1 l_2)
   \sum_{r^{(1)}+r^{(2)}-r^{(3)}-r^{(4)} = 0\atop
	 |r^{(1)}|\le |r^{(2)}|\le |r^{(3)}|\le |r^{(4)}|}
	 |\hat f_{r^{(1)}} \hat f_{r^{(2)}} 
	  \hat f_{r^{(3)}} \hat f_{r^{(4)}}| \cr
    }$$
which in turn is bounded above by
$$ \eqalignno{
   24 (l_1 l_2)
   \sum_{0\le m_1\le m_2 \le m_3} \ &
   \sum_{2^{m_1} \le |r^{(1)}| < 2^{m_1+1}} |\hat f_{r^{(1)}}|
   \sum_{2^{m_2} \le |r^{(2)}| < 2^{m_2+1}} |\hat f_{r^{(2)}}| \cr & \kern - 1cm
   \sum_{2^{m_3} \le |r^{(3)}| < 2^{m_3+1}} \ 
   I_{|r^{(1)}+r^{(2)}-r^{(3)}| \ge |r^{(3)}|} \
   |\hat f_{r^{(3)}} \hat f_{r^{(1)}+r^{(2)}-r^{(3)}}| .\cr }
   $$
In bounding this quantity, let us start by looking at the inner sum:
$$ \sum_{2^{m_3} \le |r^{(3)}| < 2^{m_3+1}}
   I_{|r^{(1)}+r^{(2)}-r^{(3)}| \ge |r^{(3)}|} \
   |\hat f_{r^{(3)}} \hat f_{r^{(1)}+r^{(2)}-r^{(3)}}| .$$
Using the Cauchy-Schwartz formula, this can be bounded above by
$$ \big(\sum_{2^{m_3} \le |r^{(3)}| < 2^{m_3+1}} 
     |\hat f_{r^{(3)}}|^2 \big)^{1/2} 
    \big(\sum_{2^{m_3} \le |r^{(3)}| < 2^{m_3+1}}
     I_{|r^{(1)}+r^{(2)}-r^{(3)}| \ge |r^{(3)}|} \
     |\hat f_{r^{(1)}+r^{(2)}-r^{(3)}}|^2
     \big)^{1/2}. $$
The first term in this product is $A_{m_3}$.  As for the second term,
since $m_1 \le m_2 \le m_3$, it follows that
$|r^{(1)}+r^{(2)}-r^{(3)}| \le 3\cdot 2^{m_3}$, and hence the second
term is bounded by $(A_{m_3}^2 + A_{m_3+1}^2 + A_{m_3+2}^2)^{1/2}$.
Thus, the above quantity can be bounded above by 
$A_{m_3}^2 + A_{m_3+1}^2 + A_{m_3+2}^2$.

Furthermore
$$ \sum_{2^{m} \le |r| < 2^{m_+1}} | \hat f_r |
   \le
   \left(\sum_{2^{m} \le |r| < 2^{m_+1}} 1 \right)^{1/2}
   \left(\sum_{2^{m} \le |r| < 2^{m_+1}} | \hat f_r |^2 \right)^{1/2}
   \le
   c 2^m A_m ,$$
since the number of points $r$ such that $2^{m} \le |r| < 2^{m+1}$ is
$2^{2m}$ to within a constant factor.

Thus
$$ \normo f_4^4 \le c
   \sum_{m_3=0}^\infty 
   (A_{m_3}^2 + A_{m_3+1}^2 + A_{m_3+2}^2)
   \sum_{m_2=0}^{m_3} 
   2^{m_2} A_{m_2}
   \sum_{m_1=0}^{m_2}
   2^{m_1} A_{m_1} .$$
Now, applying Cauchy-Schwartz, we see that the inner sum obeys the
inequalities
$$ \sum_{m_1=0}^{m_2}
   2^{m_1} A_{m_1}
   \le
   \left( \sum_{m_1=0}^{m_2} 2^{m_1} \right)^{1/2}
   \left( \sum_{m_1=0}^{m_2} 2^{m_1} A_{m_1}^2 \right)^{1/2}
   \le c 2^{m_2/2} \normo{D^{1/2} f}_2 .$$
Thus
$$ \normo f_4^4 \le c \normo{D^{1/2} f}_2
   \sum_{m_3=0}^\infty 
   (A_{m_3}^2 + A_{m_3+1}^2 + A_{m_3+2}^2)
   \sum_{m_2=0}^{m_3} 
   2^{3 m_2/3} A_{m_2} .$$
Applying Cauchy-Schwartz once more in a similar fashion, we get that
$$ \sum_{m_2=0}^{m_3} 
   2^{3 m_2/3} A_{m_2}
   \le
   c 2^{m_3} \normo{D^{1/2} f}_2 ,$$
and so
$$ \normo f_4^4 \le c \normo{D^{1/2} f}_2^2 
   \sum_{m_3=0}^\infty
   2^{m_3}
   (A_{m_3}^2 + A_{m_3+1}^2 + A_{m_3+2}^2) ,$$
from which we see that $\normo f_4^4 \le c \normo{D^{1/2} f}_2^4$ as required.

\hfill Q.E.D.

\beginsection Acknowledgments

This paper owes a lot to other people who have helped me over the years.
I would like to thank to Joel Avrin and John Gibbon
who explained the Navier-Stokes equation to me three years ago.
I would like to thank George Sell,
Mohammed Ziane and Benoit Desjardins for improvements and
corrections to this
paper.

\beginsection References

\frenchspacing

\item{A}
Avrin, Joel D. Large-eigenvalue global existence and regularity results for 
the Navier-Stokes equation. {\sl J. Differential Equations 127 (1996),
no. 2, 365--390.} 

\item{CF} Constantin, Peter; Foia\c s, Ciprian, {\sl Navier-Stokes equations.
Chicago Lectures in Mathematics. University of Chicago Press, Chicago, IL, 
1988.}

\item{DG}
Doering, Charles R.; Gibbon, J. D. {\sl Applied analysis of the Navier-Stokes 
equations. Cambridge Texts in Applied Mathematics. Cambridge
University Press, Cambridge, 1995.}

\item{I1}
Iftimie, Drago\c s Les \'equations de Navier-Stokes $3$D vues comme une 
perturbation des \'equations de Navier-Stokes $2$D. (French) [The
$3$D Navier-Stokes equations seen as a perturbation of the $2$D Navier-Stokes 
equations] {\sl C. R. Acad. Sci. Paris S\'er. I Math. 324 (1997), no. 3, 
271--274.}

\item{I2}
Iftimie, Drago\c s
The 3D Navier-Stokes equations seen as a 
perturbation of the 2D Navier-Stokes equations. {\sl preprint, available
at \hfill\break{\tt
http://www.maths.univ-rennes1.fr/\~{}iftimie/publications.html}}

\item{MTZ}
Moise, I.; Temam, R.; Ziane, M. 
Asymptotic analysis of the Navier-Stokes equations in thin domains. 
Dedicated to Olga Ladyzhenskaya. 
{\sl Topol. Methods Nonlinear Anal. {\bf 10} (1997), no. 2, 249--282.}

\item{RS1}
Raugel, Genevi\`eve; Sell, George R. Navier-Stokes equations on thin 
$3$D domains. I. Global attractors and global regularity of solutions. 
{\sl J. Amer. Math. Soc. {\bf 6} (1993), no. 3, 503--568.}

\item{RS2}
Raugel, G.; Sell, G. R. Navier-Stokes equations on thin $3$D domains. II. 
Global regularity of spatially periodic solutions. 
{\sl Nonlinear partial
differential equations and their applications. Coll\`ege de France 
Seminar, Vol. XI (Paris, 1989--1991), 205--247, Pitman Res. 
Notes Math. Ser., 299, Longman Sci.
Tech., Harlow, 1994.}

\item{RS3}
Raugel, Genevi\`eve; Sell, George R. Navier-Stokes equations in 
thin $3$D domains. III. Existence of a global attractor. 
{\sl Turbulence in fluid flows,
137--163, IMA Vol. Math. Appl., 55, Springer, New York, 1993.}

\item{S}
Stein, Elias M. {\sl Singular integrals and differentiability properties of 
functions. Princeton Mathematical Series, No. 30 Princeton University
Press, Princeton, N.J. 1970.}

\item{T}
Temam, Roger {\sl Infinite-dimensional dynamical systems in mechanics 
and physics. Second edition. Applied Mathematical Sciences, 68.
Springer-Verlag, New York, 1997.}

\item{TZ1}Temam, R.; Ziane, M. Navier-Stokes equations in three-dimensional 
thin domains with various boundary conditions. 
{\sl Adv. Differential Equations
{\bf 1} (1996), no. 4, 499--546.}

\item{TZ2}
Temam, R.; Ziane, M. 
Navier-Stokes equations in thin spherical domains. 
{\sl Optimization methods in partial differential equations (South Hadley,
MA, 1996), 281--314, Contemp. Math., {\bf 209}, Amer. Math. Soc., Providence, 
RI, 1997.}

\bigskip
\noindent 
Stephen Montgomery-Smith \hfil\break
Math.\ Dept., University of Missouri \hfil\break
Columbia, MO 65211, U.S.A. \hfil\break
email: stephen@math.missouri.edu \hfil\break
http://math.missouri.edu/$\sim$stephen

\bigskip

\beginsection Addendum:  July 8 1999.

On page~3, in both Theorems~1 and~2, the sentence that begins ``{\sl
If $\displaystyle t \ge c \, {l_1^2 \over \nu}$}~''
should read

{\narrower \smallskip \noindent \sl 
Furthermore
$$ \limsup_{t \to \infty }\normo{u(t)}_{H^1} 
  \le c \max\left\{{l_1 \over \nu} F,
   {l_1^{7/2} \over \nu^3 l_2^{1/2} } \, \epsilon^{-1/2} F^2 \right\}. $$
\smallskip}

Similarly, for Lemma's~3 and~5 (pages~8/9 and 14 respectively), the
phrase that begins ``{\sl and if $t \ge c_{16}$ then}'' and ends
``{\sl for $0 \le t < T$}'' (respectively ``{\sl for $0 \le t < \infty$}'') 
should be replaced by

{\narrower \smallskip \noindent \sl
and
$$ \limsup_{t \to \infty} \psi \le c_{17} 
      \max\{\epsilon^{-1/2} F^2, F \}. \eqno (4.12) $$
\smallskip}
\noindent For Lemma~5 the corresponding equation number is (8.12).

The second to last paragraph of the proof of Lemma~3 (page~9) should
be changed to

{\narrower\smallskip\noindent
To obtain (4.12) is similar.  We see that given $\epsilon>0$ there are positive
numbers
$\tau$, $c_{23}$ and $c_{24}$ (where $\tau$ depends upon $\epsilon$)
such that if $t \ge \tau$ then
$\phi \le c_{23} F+\epsilon$ and $\theta \le c_{24} F+\epsilon$.  
Apply the above argument,
except integrate from $\tau$ to $t$ instead of from $0$ to $t$.\smallskip
}

\bye